\documentclass[openacc]{rsproca_new}

\usepackage{natbib}
\usepackage{xcolor}
\usepackage{tikz}
\usepackage{pgfplots}
\usepgfplotslibrary{colormaps}
\usetikzlibrary{shapes}
\usetikzlibrary{pgfplots.colormaps}
\newcommand\mtcaption{\fontsize{8pt}{8pt}\selectfont}
\newcommand\mtlabel{\fontsize{8pt}{8pt}\selectfont}
\pgfplotsset{compat=1.16}
\pgfplotsset{every tick label/.append style={font=\mtlabel},
label style={font=\mtcaption,}}

\usetikzlibrary{arrows,matrix,positioning,fit}

\tikzstyle{process} = [rectangle, draw, text centered, minimum height=2em]
\tikzstyle{decision} = [diamond, draw, text centered, minimum height=2em, aspect=3]
\tikzstyle{connector} = [draw, -latex']

\newcommand{\newtext}[1]{#1}
\newcommand{\nphy}{k}
\newcommand{\nnum}{n}
\newcommand{\nint}{m}

\begin{document}

\title{Real-time hybrid testing using iterative control for periodic oscillations}

\author{
S. Beregi$^{1,2}$, D.A.W. Barton$^{1}$, D. Rezgui$^{1}$ and S.A. Neild$^{1}$}

\address{$^{1}$Faculty of Engineering, University of Bristol, UK\\
$^{2}$Department of Infectious Disease Epidemiology, Imperial College London, UK\\}

\subject{xxxxx, xxxxx, xxxx}


\corres{S. Beregi\\
\email{sberegi@ic.ac.uk}}

\begin{abstract}
Real-time hybrid testing is a method in which a substructure of the system is realised experimentally and the rest numerically. The two parts interact in real time to emulate the dynamics of the full system. Such experiments however are often difficult to realise as the actuators and sensors, needed to ensure compatibility and force-equilibrium conditions at the interface, can seriously affect the predicted dynamics of the system and result in stability and fidelity issues. The traditional approach of using feedback control to overcome the additional unwanted dynamics is challenging due to the presence of an outer feedback loop, passing interface displacements or forces to the numerical substructure. We, therefore, advocate for an alternative approach, removing the problematic interface dynamics with an iterative scheme to minimise interface errors, thus, capturing the response of the true assembly. The technique is examined by hybrid testing of a bench-top four-storey building with different interface configurations, where using conventional hybrid measurement techniques is very challenging. A case where the physical part exhibits nonlinear restoring force characteristics is also considered. These tests show that the iterative approach is capable of handling even scenarios which are theoretically infeasible with feedback control.
\end{abstract}


\begin{fmtext}

\end{fmtext}
\maketitle
\section{Introduction}

Carrying out experiments on large-scale structures is challenging in engineering, as the instrumentation, production, moving and transportation of large prototypes is usually very costly and time-consuming. Thus, in many cases, a full-scale physical experiment it is not feasible. To address these challenges, one may conduct a physical test on particular substructures of a large-scale assembly. However, the tested substructure is likely to behave differently to the full-scale tests because boundary conditions are not consistent with the full structure.

In contrast, numerical simulations of engineering structures are inexpensive and there is much flexibility to try different concepts and designs for a product. While a full-numerical test is usually cost-effective, physical experiments are required for validating these `numerical twins' of the physical structures. Thus, it is not entirely possible to eliminate the challenges with physical experiments mentioned above. 

Real-time hybrid testing (a {sub-field within} dynamic substructuring; \cite{subsAllen}) is a part experimental, part numerical approach to obtain the dynamic response of an engineering structure, where only a part of the tested assembly is physically realised, whereas the other part is simulated. It naturally follows from this concept that, in order to replicate the dynamic behaviour of the true assembly, the physical substructure has to interact in real-time with the numerical half of the model.

A hybrid experiment may be useful in several scenarios, especially during the design of an engineering product. A physical-numerical hybrid test can be a more flexible alternative to a full-scale experiment since it allows for testing components in different environments without the need to physically realise every scenario or environment considered of interest. Thus, the hybrid testing concept can enable significant savings in experimental costs.

Thanks to its advantages, hybrid testing sparked significant interest in engineering and several previous experimental studies were carried out to investigate dynamic substructuring in various scenarios. The range of use cases includes civil engineering (\cite{Qian_Dyke,Gao_Dyke,Miraglia}), mechanical- (\cite{vdSeijs, Nicgorski, Meggitt}) or vibro-acoustic systems (\cite{RixenGuitar}). There are also examples of use with problems including coupled physics, e.g. structural and aero- or fluid dynamics, such as in case of an aircraft wing \cite{Ruffini2020} or floating wind turbines \cite{Bachynski}.

Real-time hybrid testing takes a similar approach to hardware-in-the-loop experiments, where a physical structure interacts with a simulated environment. While the two concepts are similar, the contrast between the two techniques can be found in the areas of application. Typical examples for the hardware-in-the-loop testing approach are mechatronic (electrical-mechanical) system testing, e.g. flight or driving simulators \cite{HIL}. In principle, adopting the same idea to test mechanical structures in a hybrid experimental environment is straightforward. However, the interface dynamics in mechanical hybrid experiments is often non-negligible and has a significant effect on the overall test dynamics \cite{Horiuchi}.

\begin{figure}[!h]
    \centering\includegraphics[]{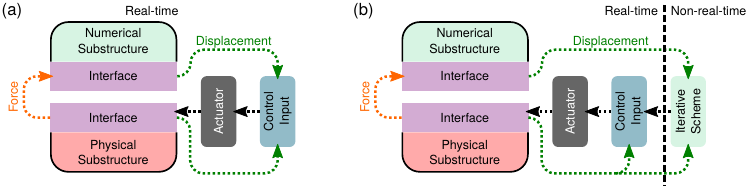}
    \caption{The real-time hybrid testing concept. Panel (a): Substructuring with a full real-time feedback-loop. Panel (b): Substructuring with real-time force-feedback and displacement control with non-real-time iteration. The force measured at the interface of the physical substructure is used within the numerical substructure to calculate the required displacement at the interface.}
    \label{fig:substr}
\end{figure}

The advantages of real-time hybrid testing though come with challenges as well, which make the realisation of such hybrid experiments difficult. In mechanical structures, real-time data from the physical substructure is typically collected by electrodynamic transducers, e.g. piezoelectric load cells or accelerometers, and fed to the simulated part of the system. In the opposite direction across the interface, actuators are used to impose the behaviour of the numerical substructure on the physical experiment. This leads to a feedback-loop as illustrated by panel (a) in Fig.~\ref{fig:substr}. Note that, here, we depict a displacement control based set-up, it is also possible to feed displacement to the numerical substructure and then use the actuator to impose a force on the physical substructure. 

Interface errors inevitably lead to fidelity issues, since the behaviour of the true hybrid system with the interface dynamics is always different (to some degree) than that of the true assembly. In \cite{Insam2022}, the fidelity issue in real-time substructuring tests is considered and an error quantification method is proposed. Nevertheless, in some cases, one may not just simply perceive an error in the response of the hybrid system, but the contribution of the actuator dynamics often leads to instabilities related to the feedback-loop which can make the hybrid experiment potentially infeasible. Instabilities in substructuring may be encountered even in conceptually simple systems, cf., previous studies on theoretical substructurability of low-degree-of-freedom mass-spring-damper systems \cite{Terkovics2016}, a coupled oscillator-pendulum system \cite{Tu}, and cantilever-beam with a PDE model \cite{KyrychkoPDE}. In analyses focusing on transducer-induced instabilities, it is commonplace to approximate the interface error as a time-delayed response and perform stability analysis on the resulting system of delay differential equations \cite{HaleLunel}.

The standard way of reducing interface error is using interface control \cite{Ruffini2020,Gawthrop2007} and/or an inverse-model-based compensation technique, relying on the thorough identification of the actuator dynamics. By means of this approach, one may employ model predictive control (see \cite{forwardpred} and \cite{Tsokanas2022}) to compensate the phase lag induced by the actuators. For time-delayed systems in particular, one may also consider using the technique of finite spectrum assignment \cite{FSA}. This method however can be very cumbersome, as it requires thorough system identification, which can involve inaccuracies and has to be performed repeatedly every time the experimental setup is changed.

To tackle the above limitations, we propose an alternative, iterative approach for hybrid experiments excited periodically. This approach removes the interface from the feedback-loop, thus, avoiding potential stability issues involved. Instead, we provide the controller with a periodic demand and monitor the response, updating the harmonic amplitudes of the demand in non-real-time until the interface compatibility and equilibrium is ensured. {The closest existing work in the literature is} \cite{Witteveen}, where they successfully realise a non-simultaneous hybrid experiment, testing a physical shock absorber with a numerical wheel-suspension model. In their study, Witteveen et al. transform the coupling errors into the frequency domain, and use an iterative technique to establish the required forcing input to eliminate it. Note that the concept of hybrid testing in the frequency domain is also well-established in the literature {with off-line FRF-based approaches such as} \cite{RixenImpulse,Geradin} and,  modal approaches \cite{AllenISMA}; it is one of the standard techniques in the field.

Here we propose a different approach inspired by control-based continuation \cite{Sieber2008}, a technique for experimentally tracking limit cycles in nonlinear systems (both stable and unstable). The method employs a iteration to find the control target ensuring that the control converges to the same steady-state response as would be the case in the original system.

By capturing the FRF of the true assembly for an excitation with a stepped frequency, we treat the hybrid experiment as a varying parameter problem. We demonstrate the usefulness of this feature, by taking advantage of the original idea of the iterative approach coming from experimental bifurcation analysis. Thus, the method is capable of performing hybrid testing for structures with a nonlinear restoring force characteristics.

In this study, we use the example of a bench-top-sized 4-storey building to demonstrate the versatility and effectiveness of the iterative hybrid testing technique. Similar structures have been already used in studies on dynamic substructuring (e.g. \cite{Tian}) exploiting the easily describable dynamic response, namely, the behaviour of the 4-storey structure can be well-characterised by a 4 degree-of-freedom (DoF) model with concentrated mass, damping and stiffness parameters \cite{DalBorgo3sorey,GardnerSHM,GardnerBayes}. The main challenge is that imperfect control compensation for the actuator dynamics leads to a mismatch across the interface and hence a misrepresentation of the dynamics of the structure being emulated. We first show that even for this conceptually simple structure, real-time feedback control of the hybrid test is very challenging and infeasible for some parameter ranges. We then demonstrate the effectiveness and versatility of the new proposed approach by performing dynamic substructuring experiments. 

The rest of the paper is constructed as follows. In Section 2, we present the mathematical model of real-time hybrid testing with time-delayed-feedback from the numerical substructure and analyse the theoretical substructurability by performing a stability analysis on the hybrid system. In Section 3, we introduce the iterative algorithm eliminating direct feedback and the related instability from the structure. Section 4 presents the experimental setup which was used to produce the experimental results shown in Section 5. We first demonstrate successful implementation of the iterative hybrid testing approach on the linear system before extending the study to a weakly nonlinear version of the structure. The conclusions are given in Section 6.






\section{Mathematical modelling}

\subsection{The time-delayed model of a real-time hybrid test}\label{ch:substr}

In this section, a mathematical model of the full hybrid system will be developed and examined to explore the control challenge discussed in the introduction. We note that for a real test, a full model of the system may not be available. \newtext{The underlying purpose of this section is to demonstrate that the type of hybrid system explored here exhibits a particularly challenging type of instability, namely a \emph{neutral delay}, whereby there are parameter regimes in which the system is unconditionally unstable --- any delay renders the system unstable. This provides motivation for the iterative control scheme presented here, which circumvents issues of delay.}

The motivation to conduct a hybrid test is often due to challenges with modelling components of the physical substructure, hence it being represented experimentally. Here, to allow us to analyse the hybrid test approach, we will consider a system that can be represented numerically. 

We consider an $l$-degree-of-freedom (DoF) linear dynamical system as the \emph{true assembly} we would like to emulate with a part physical, part numerical experiment. Note that hybrid testing may be realised for more complex, in particular, nonlinear models; nevertheless, we consider a general linear case as it allows us to highlight the challenge without additional algebraic complexity.

The position and orientation of the true assembly is described by the vector of generalised coordinates $\mathbf{q}(t) \in \mathbb{R}^l$ where the individual generalised coordinates $q_i(t)$ are usually chosen to fully define the system with the minimal number of coordinates. Thus, the equations of motion can be given in the matrix form
\begin{equation}\label{eq:true_assembly}
    \mathbf{M}\ddot{\mathbf{q}}(t)+\mathbf{C}\dot{\mathbf{q}}(t)+\mathbf{K}\mathbf{q}(t) = \mathbf{F}(t),
\end{equation}
where $\mathbf{M,C}$ and $\mathbf{K}$ are the mass, damping and stiffness matrices, respectively, ($\mathbf{M,C,K} \in \mathbb{R}^{l \times l}$), whereas $\mathbf{F}(t) \in \mathbb{R}^l$ describes external forcing applied to the structure.

In the hybrid experiment, only a part of the true assembly is represented in the \emph{physical substructure} -- this does not require a model in the hybrid test but to analyse the system, it will be represented by $\nphy$ of the DoFs in the full model -- the rest is replaced by a numerical model, referred to as the $\nnum$ DoF \emph{numerical substructure}. The compatibility of the physical and the numerical substructures is ensured by $\nint$ shared degrees-of-freedoms. The corresponding rigid bodies are referred to as the \emph{interface} and distinguished from the \emph{internal} parts of the two substructures. As such, the sum of the DoF of the two substructures is larger than $l$ (the DoF of the true assembly) by $\nint$, i.e. $l=\nphy+\nnum-\nint$.

The hybrid testing concept assumes that the two parts are clearly separable which is ensured by the coefficient matrices $\mathbf{M,C}$ and $\mathbf{K}$ in Eq.~\eqref{eq:true_assembly} having a block diagonal structure

\begin{equation}
    \begin{tikzpicture}[baseline=(current  bounding  box.center)]
        \matrix [matrix of math nodes,left delimiter=(,right delimiter=)] (m)
        {
            {\bullet}_{11} & {\bullet}_{12} & {\bf 0}\\
            {\bullet}_{21} & {\bullet}_{22} & {\bullet}_{23}\\
            {\bf 0} & {\bullet}_{32} & {\bullet}_{33} \\
        };
        \node[draw,red,dashed,fit=(m-1-1) (m-2-2),inner sep=0pt]{};
        \node[draw,green,dashed,fit=(m-2-2) (m-3-3),inner sep=0pt]{};

    \end{tikzpicture},
\end{equation}
where $\bullet_{11} \in \mathbb{R}^{(\nphy-\nint) \times (\nphy-\nint)}$, $\bullet_{22} \in \mathbb{R}^{\nint \times \nint}$ and $\bullet_{33} \in \mathbb{R}^{(\nnum-\nint) \times (\nnum-\nint)}$. In the formulae presented here, the bullets stand for the corresponding block of $\mathbf{M,C}$ and $\mathbf{K}$; e.g., for the mass matrix, $\bullet_{11}$ is interpreted as $\mathbf{M}_{11}$. We mark the matrix-blocks corresponding to the physical and numerical substructures with red, green boxes, respectively.
Notice that the interface element $\bullet_{22}$ of the coefficient matrices is split between the physical and the numerical substructures. Thus, the two substructures should be considered such that the sum of the corresponding elements is equal to that of the true assembly
$\bullet_{22} = \bullet_{22}^{\rm phy} + \bullet_{22}^{\rm num}$.

The equations of motion of the physical and numerical substructures can be therefore expressed as
\begin{equation}\label{eq:phy}
    \begin{tikzpicture}[baseline=(current  bounding  box.center),
        every left delimiter/.style={xshift=1.0em},
        every right delimiter/.style={xshift=-1.0em},]
        \matrix [matrix of math nodes,left delimiter=(,right delimiter=)] (m)
        {
            {\bf{M}}_{11} & {\bf{M}}_{12}\\
            {\bf{M}}_{21} & {\bf{M}}_{22}^{\rm phy}\\
        };
    \end{tikzpicture}
    \begin{tikzpicture}[baseline=(current  bounding  box.center),
        every left delimiter/.style={xshift=1.0em},
        every right delimiter/.style={xshift=-1.0em},]
        \matrix [matrix of math nodes,left delimiter=(,right delimiter=)] (m)
        {
            {\ddot{\bf{x}}}_{\rm phy}(t)\\
            {\ddot{\bf{x}}}_{\rm int}(t)\\
        };
    \end{tikzpicture}
    +
    \begin{tikzpicture}[baseline=(current  bounding  box.center),
        every left delimiter/.style={xshift=1.0em},
        every right delimiter/.style={xshift=-1.0em},]
        \matrix [matrix of math nodes,left delimiter=(,right delimiter=)] (m)
        {
            {\bf{C}}_{11} & {\bf{C}}_{12}\\
            {\bf{C}}_{21} & {\bf{C}}_{22}^{\rm phy}\\
        };
    \end{tikzpicture}
    \begin{tikzpicture}[baseline=(current  bounding  box.center),
        every left delimiter/.style={xshift=1.0em},
        every right delimiter/.style={xshift=-1.0em},]
        \matrix [matrix of math nodes,left delimiter=(,right delimiter=)] (m)
        {
            {\dot{\bf{x}}}_{\rm phy}(t)\\
            {\dot{\bf{x}}}_{\rm int}(t)\\
        };
    \end{tikzpicture}
    +
    \begin{tikzpicture}[baseline=(current  bounding  box.center),
        every left delimiter/.style={xshift=1.0em},
        every right delimiter/.style={xshift=-1.0em},]
        \matrix [matrix of math nodes,left delimiter=(,right delimiter=)] (m)
        {
            {\bf{K}}_{11} & {\bf{K}}_{12}\\
            {\bf{K}}_{21} & {\bf{K}}_{22}^{\rm phy}\\
        };
    \end{tikzpicture}
    \begin{tikzpicture}[baseline=(current  bounding  box.center),
        every left delimiter/.style={xshift=.7em},
        every right delimiter/.style={xshift=-.7em},]
        \matrix [matrix of math nodes,left delimiter=(,right delimiter=)] (m)
        {
            {{\bf{x}}}_{\rm phy}(t)\\
            {{\bf{x}}}_{\rm int}(t)\\
        };
    \end{tikzpicture}
    =
    \begin{tikzpicture}[baseline=(current  bounding  box.center),
        every left delimiter/.style={xshift=.7em},
        every right delimiter/.style={xshift=-.7em},]
        \matrix [matrix of math nodes,left delimiter=(,right delimiter=)] (m)
        {
            {{\bf{F}}}_{\rm phy}(t)\\
            {{\bf{F}}}_{\rm int}^{\rm phy}(t)\\
        };
    \end{tikzpicture},
\end{equation}
and 
\begin{equation}\label{eq:virt}
    \begin{tikzpicture}[baseline=(current  bounding  box.center),
        every left delimiter/.style={xshift=1.0em},
        every right delimiter/.style={xshift=-1.0em},]
        \matrix [matrix of math nodes,left delimiter=(,right delimiter=)] (m)
        {
            {\bf{M}}_{22}^{\rm num} & {\bf{M}}_{32}\\
            {\bf{M}}_{23} & {\bf{M}}_{33}\\
        };
    \end{tikzpicture}
    \begin{tikzpicture}[baseline=(current  bounding  box.center),
        every left delimiter/.style={xshift=1.0em},
        every right delimiter/.style={xshift=-1.0em},]
        \matrix [matrix of math nodes,left delimiter=(,right delimiter=)] (m)
        {
            {\ddot{\pmb{\xi}}}_{\rm int}(t)\\
            {\ddot{\pmb{\xi}}}_{\rm num}(t)\\
        };
    \end{tikzpicture}
    +
    \begin{tikzpicture}[baseline=(current  bounding  box.center),
        every left delimiter/.style={xshift=1.0em},
        every right delimiter/.style={xshift=-1.0em},]
        \matrix [matrix of math nodes,left delimiter=(,right delimiter=)] (m)
        {
            {\bf{C}}_{22}^{\rm num} & {\bf{C}}_{23}\\
            {\bf{C}}_{32} & {\bf{C}}_{33}\\
        };
    \end{tikzpicture}
    \begin{tikzpicture}[baseline=(current  bounding  box.center),
        every left delimiter/.style={xshift=1.0em},
        every right delimiter/.style={xshift=-1.0em},]
        \matrix [matrix of math nodes,left delimiter=(,right delimiter=)] (m)
        {
            {\dot{\pmb{\xi}}}_{\rm int}(t)\\
            {\dot{\pmb{\xi}}}_{\rm num}(t)\\
        };
    \end{tikzpicture}
    +
    \begin{tikzpicture}[baseline=(current  bounding  box.center),
        every left delimiter/.style={xshift=1.0em},
        every right delimiter/.style={xshift=-1.0em},]
        \matrix [matrix of math nodes,left delimiter=(,right delimiter=)] (m)
        {
            {\bf{K}}_{22}^{\rm num} & {\bf{K}}_{23}\\
            {\bf{K}}_{32} & {\bf{K}}_{33}\\
        };
    \end{tikzpicture}
    \begin{tikzpicture}[baseline=(current  bounding  box.center),
        every left delimiter/.style={xshift=.7em},
        every right delimiter/.style={xshift=-.7em},]
        \matrix [matrix of math nodes,left delimiter=(,right delimiter=)] (m)
        {
            {\pmb{\xi}}_{\rm int}(t)\\
            {\pmb{\xi}}_{\rm num}(t)\\
        };
    \end{tikzpicture}
    =
    \begin{tikzpicture}[baseline=(current  bounding  box.center),
        every left delimiter/.style={xshift=.7em},
        every right delimiter/.style={xshift=-.7em},]
        \matrix [matrix of math nodes,left delimiter=(,right delimiter=)] (m)
        {
            {{\bf{F}}}_{\rm int}^{\rm num}(t)\\
            {{\bf{F}}}_{\rm num}(t)\\
        };
    \end{tikzpicture},
\end{equation}
respectively. In the equations above, $\mathbf{x}_{\rm phi}$ and $\mathbf{x}_{\rm int}$ contain the generalised coordinates corresponding to the internal part and the interface of the physical substructure,
\begin{equation}
    \mathbf{x}_{\rm phy} = \left( \begin{matrix} x_1 & ... & x_{\nphy-\nint} \end{matrix} \right)^{\intercal},
    \quad
    \mathbf{x}_{\rm int} = \left( \begin{matrix} x_{\nphy-\nint+1} & ... & x_\nphy \end{matrix} \right)^{\intercal},
\end{equation}
and similarly, $\pmb{\xi}_{\rm int}$ and $\pmb{\xi}_{\rm num}$ represent the generalised coordinates of the interface and the internal part of the numerical substructure
\begin{equation}
    \pmb{\xi}_{\rm int} = \left( \begin{matrix} \xi_{\nphy-\nint+1} & ... & \xi_\nphy \end{matrix} \right)^{\intercal},
    \quad
    \pmb{\xi}_{\rm num} = \left( \begin{matrix} \xi_{\nphy+1} & ... & \xi_{l} \end{matrix} \right)^{\intercal}.
\end{equation}

To replicate the behaviour of the true assembly by the hybrid system, the motion of the interface in the physical and numerical substructures should be synchronised. In the substructurability analysis we consider the approach of feedback control, using force-feedback from the physical to the numerical substructure, and displacement feedback from the numerical to the physical substructure, cf. panel (a) in Fig.~\ref{fig:substr}. Since sensors tend to introduce negligible phase-lag, we assume the force-feedback from the physical to the numerical substructure to be real-time; thus, the corresponding physical and numerical coupling forces are equal and point in the opposite direction
\begin{equation}\label{eq:cpl_forc}
    \mathbf{F}_{\rm int}^{\rm phy}(t) = -\mathbf{F}_{\rm int}^{\rm num}(t)
\end{equation}
Conversely, we consider a time-delayed displacement feedback from the numerical to the physical subsystem, which can be expressed by the delayed displacement coupling condition
\begin{equation}\label{eq:cpl_disp}
    \mathbf{x}_{\rm int}(t) = \pmb{\xi}_{\rm int}(t-\tau),
\end{equation}
where $\tau$ denotes the feedback time-delay. The physical-background of the delayed feedback is rooted in sensor-dynamics, sampling and signal processing times, filtering delays, and actuator dynamics. While the contribution of these effects can differ in every experiment, since the aim of our study is to demonstrate the challenges of feedback-based real-time hybrid testing in the qualitative sense, we consider a single time delay in our model.

Expressing the coupling force from the second row of Eq.~\eqref{eq:phy}, and substituting in the delayed displacement coupling condition we obtain
\begin{multline}
    \mathbf{F}_{\rm int}^{\rm phy}  = 
    \mathbf{M}_{21} \ddot{\mathbf{x}}_{\rm phy}(t)+\mathbf{M}_{22}^{\rm phy} \ddot{\pmb{\xi}}_{\rm int}(t-\tau) \\
    + \mathbf{C}_{21} \dot{\mathbf{x}}_{\rm phy}(t)+\mathbf{C}_{22}^{\rm phy} \dot{\pmb{\xi}}_{\rm int}(t-\tau) +
    \mathbf{K}_{21} {\mathbf{x}}_{\rm phy}(t)+\mathbf{K}_{22}^{\rm phy} {\pmb{\xi}}_{\rm int}(t-\tau).
\end{multline}
Then, using the force coupling condition, we can substitute the expression above into the right-hand-side of Eq.~\eqref{eq:virt}. As a result, using the first row of Eq.~\eqref{eq:phy} and the whole matrix formula in Eq.~\eqref{eq:virt} we can compose the equation of motion for substructuring with delayed displacement-feedback
\begin{equation}\label{eq:eom_substr}
    \mathbf{M}_{\tau} \ddot{\mathbf{q}}(t-\tau) +
    \mathbf{M}_{0} \ddot{\mathbf{q}}(t)+
    \mathbf{C}_{\tau} \dot{\mathbf{q}}(t-\tau) +
    \mathbf{C}_{0} \dot{\mathbf{q}}(t)+
    \mathbf{K}_{\tau} {\mathbf{q}}(t-\tau) +
    \mathbf{K}_{0} {\mathbf{q}}(t) = \mathbf{F}_{\rm substr}(t),
\end{equation}
where, the vector of generalised coordinates is given by
\begin{equation}
 \mathbf{q}(t) = \left( \begin{matrix} \mathbf{x}_{\rm phy}(t) & \pmb{\xi}_{\rm int}(t) & \pmb{\xi}_{\rm num}(t) \end{matrix} \right)^{\intercal},
\end{equation}
while using bullets to represent the corresponding non-zero elements of each, the delayed and real-time mass, damping and stiffness matrices $\mathbf{M}_{\tau},\mathbf{C}_{\tau},\mathbf{K}_{\tau},\mathbf{M}_{0},\mathbf{C}_{0}$ and $\mathbf{K}_{0}$ can be given in the structure
\begin{equation}
    \bullet_{\tau} =
    \begin{tikzpicture}[baseline=(current  bounding  box.center),
        every left delimiter/.style={xshift=1.0em},
        every right delimiter/.style={xshift=-1.0em},]
        \matrix [matrix of math nodes,left delimiter=(,right delimiter=)] (m)
        {
            {\bf{0}} & {\bullet}_{12} & {\bf{0}}\\
            {\bf{0}} & {\bullet}_{22}^{\rm phy} & {\bf{0}}\\
            {\bf{0}} & {\bf{0}} & {\bf{0}} \\
        };
    \end{tikzpicture},
    \quad
    \bullet_{0} = 
    \begin{tikzpicture}[baseline=(current  bounding  box.center),
        every left delimiter/.style={xshift=1.0em},
        every right delimiter/.style={xshift=-1.0em},]
        \matrix [matrix of math nodes,left delimiter=(,right delimiter=)] (m)
        {
            {\bullet}_{11} & {\bf{0}} & {\bf{0}}\\
            {\bullet}_{22} & {\bullet}_{22}^{\rm num} & {\bullet}_{23}\\
            {\bf{0}} & {\bullet}_{32} & {\bullet}_{33} \\
        };
    \end{tikzpicture},
\end{equation}
while the forcing vector reads
\begin{equation}
    \mathbf{F}_{\rm substr}(t) = \left( \begin{matrix} \mathbf{F}_{\rm phy}(t) \\ \pmb{0} \\ \mathbf{F}_{\rm num}(t)\\ \end{matrix} \right).
\end{equation}

As a result of the time-delayed displacement feedback, time delays appears in the states and time derivatives of state variables both in the equation of motion Eq.~\eqref{eq:eom_substr}, leading to differential equations of neutral type.
In a neutral type equation, the infinitely many eigenvalues accumulate at a constant ${\Re}(\lambda_i)$ in the complex plane as $i \rightarrow \infty$ (\cite{KyrychkoNDDE}). 
This means that stable parameter-combinations are only possible if the asymptote at which the eigenvalues accumulate is in the negative-real half-plane. This, as we will demonstrate, gives rise to a condition on the mass matrices $\mathbf{M}_0$ and $\mathbf{M}_{\tau}. $ As such, it is expected that only certain mass-ratios will be substructurable, whereas the system will be unconditionally unstable for others.

\subsection{Theoretical substructurability analysis of multi-storey buildings}

\begin{figure}[!h]
    \centering\includegraphics{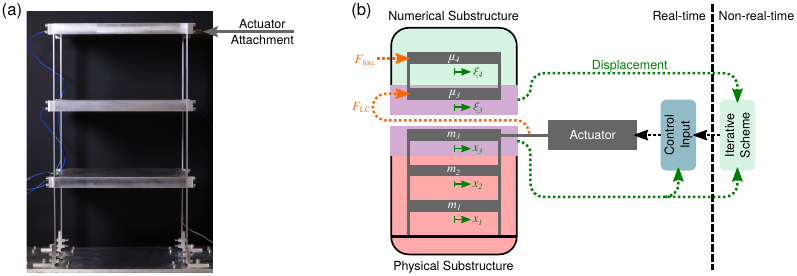}
    \caption{Real-time hybrid testing setup of a bench-top-sized four-storey building. The physical substructure is a three-storey building coupled to the numerical model of the fourth storey. The third storey serves as the interface between the substructures. Displacement control on the physical substructure is realised by an electrodynamic shaker accounting for the largest part of the time-delay in the mathematical model.}
    \label{fig:substr_with_exp}
\end{figure}

Our demonstrative example for real-time hybrid testing is a bench-top-sized four-storey building model. For the purpose of the substructurability analysis with delayed feedback, we consider a 4 degree-of-freedom lumped parameter model of the structure. The true assembly is split into two substructures at the third storey. The first and second storeys correspond to the physical substructure, the fourth storey to the numerical substructure, whereas the third storey serves as the interface. The model of the substructured assembly is depicted in Fig.~\ref{fig:substr_with_exp}. Following the notation introduced in Section 2.(a), this means that out of the $l=4$ DoF of the true assembly, $\nphy = 3$ is associated with the physical part, $\nnum = 2$ with the numerical part, while third storey's displacement being an interface coordinate ($\nint = 1$).

In the equation of motion, we refer by $m_r,\;c_r$ and $k_r$, to the masses, viscous damping and stiffness parameters of the physical side of the model, with the indices $r=1,2,3$ pointing to the storey each concentrated parameter element corresponds to. The equivalent parameters on the numerical side are denoted by $\mu_r$, $\gamma_r$ and $\sigma_r$ ($r=3,4$), for mass, damping and stiffness, respectively. The displacements of the physical storeys are given by the generalised coordinates $x_r$, while the numerical displacements are referred to as $\xi_r$.

To describe the motion of the substructured system, the following generalised coordinates are used: the physical displacements $x_1(t)$ and $x_2(t)$ and the numerical displacements $\xi_3$ and $\xi_4$. These are represented in the vector of generalised coordinates.

\begin{equation}
    \mathbf{q}(t) = \left( \begin{matrix} x_1(t) \\
    x_2(t) \\ \xi_3(t) \\ \xi_4(t) \end{matrix} \right)
\end{equation}

We assume that the displacement $x_3(t)$ of the third storey of the physical structure is constrained to the displacement of its numerical counterpart.

Thus, the equations of motion of the true assembly can be given in the matrix form as seen in Eq.~\eqref{eq:true_assembly}.
The mass matrix reads
\begin{equation}
    \label{eq:mass}
    \mathbf{M} =
    \begin{tikzpicture}[baseline=(current  bounding  box.center)]
        \matrix [matrix of math nodes,left delimiter=(,right delimiter=)] (m)
        {
            m_1 & 0 & 0 & 0\\
            0 & m_2 & 0 & 0\\
            0 & 0 & {\color{red} m_3}+{\color{green}\mu_3} & 0 \\
            0 & 0 & 0 & \mu_4 \\
        };
        \node[draw,red,dashed,fit=(m-1-1) (m-3-3),inner sep=0pt]{};
        \node[draw,green,dashed,fit=(m-3-3) (m-4-4),inner sep=0pt]{};
        \node[draw,magenta,fit=(m-3-3) (m-3-3),inner sep=0pt]{};

    \end{tikzpicture},
\end{equation}
whereas the damping and stiffness matrices are given by
\begin{equation}
    \label{eq:damping+stiffness}
    \mathbf{C} =
    \begin{tikzpicture}[baseline=(current  bounding  box.center)]
        \matrix [matrix of math nodes,left delimiter=(,right delimiter=)] (m)
        {
            c_1+c_2 & -c_2 & 0 & 0\\
            -c_2 & c_2+c_3 & -c_3 & 0\\
            0 & -c_3 & {\color{red} c_3}+{\color{green}\gamma_4} & -\gamma_4 \\
            0 & 0 & -\gamma_4 & \gamma_4 \\
        };
        \node[draw,red,dashed,fit=(m-1-1) (m-3-3),inner sep=0pt]{};
        \node[draw,green,dashed,fit=(m-3-3) (m-4-4),inner sep=0pt]{};
        \node[draw,magenta,fit=(m-3-3) (m-3-3),inner sep=0pt]{};

    \end{tikzpicture},
    \quad
    \mathbf{K} =
    \begin{tikzpicture}[baseline=(current  bounding  box.center)]
        \matrix [matrix of math nodes,left delimiter=(,right delimiter=)] (m)
        {
            k_1+k_2 & -k_2 & 0 & 0\\
            -k_2 & k_2+k_3 & -k_3 & 0\\
            0 & -k_3 & {\color{red} k_3}+{\color{green}\sigma_4} & -\sigma_4 \\
            0 & 0 & -\sigma_4 & \sigma_4 \\
        };
        \node[draw,red,dashed,fit=(m-1-1) (m-3-3),inner sep=0pt]{};
        \node[draw,green,dashed,fit=(m-3-3) (m-4-4),inner sep=0pt]{};
        \node[draw,magenta,fit=(m-3-3) (m-3-3),inner sep=0pt]{};

    \end{tikzpicture}.
\end{equation}
We mark the parts of the interface term belonging to the physical and numerical substructures by red and green, respectively.
Using the derivations presented in Section \ref{ch:substr}, one can give the governing equations of the substructured system in the same form as indicated by Eq.~\eqref{eq:eom_substr}. 

\subsection{Stability analysis}

\begin{figure}[!h]
    \centering\includegraphics[width=\textwidth]{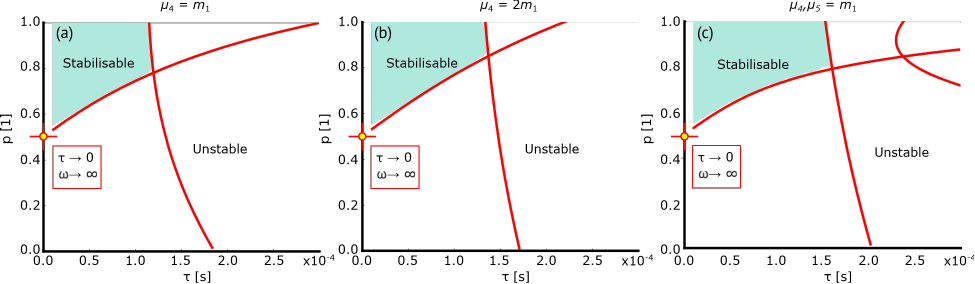}
    \caption{Stabilisability diagrams of three real-time hybrid tests. Panels (a) and (b) use a four-storey structure, while panel (c) uses a five-storey structure {with the same physical substructure as the four-storey structure and an extra storey in the numerical substructure. All tests use an interface at the third storey. The five-storey model is included to illustrate that the four-storey model is not a special case and the existence of the identified unstable regions is generic. } Stabilisable combinations of the time-delay $\tau$ and the mass-ratio $p = \mu_3 /m_3$ are shaded in green. The red curves indicate boundaries where a pair of complex conjugate characteristic exponents cross the imaginary axis. The yellow marker indicates the theoretical limit boundary of stabilisability $p=1/2$ for $\tau \rightarrow 0, \omega \rightarrow \infty$ \newtext{(see Appendix \ref{sec:asymptotics})}. We draw attention to the scale of the horizontal axis being $10^{-4}$; these conditions are extremely challenging to achieve experimentally.}
    \label{fig:stabcharts}
\end{figure}

The substructurability of the four-storey building is assessed by the stability analysis of the trivial equilibrium $\mathbf{q}(t) \equiv \mathbf{0}$ of the system given in \eqref{eq:eom_substr}. By using the exponential trial solution $\mathbf{q}(t) = \mathbf{c} e^{\lambda t}$ \cite{HaleLunel}, one obtains the characteristic equation $D(\lambda) = 0$, where the characteristic function reads
\begin{equation}\label{eq:char_fun}
    D(\lambda) = {\rm det} \left( \lambda^2 \mathbf{M}_{\tau} e^{-\lambda \tau} +\lambda^2 \mathbf{M}_0 + \lambda \mathbf{C}_{\tau} e^{-\lambda \tau} +\lambda \mathbf{C}_0 + \mathbf{K}_{\tau}e^{-\lambda \tau}+\mathbf{K}_0 \right).
\end{equation}

Stability boundaries may occur where the system has one or more characteristic exponent(s) with a zero real part, i.e. $\Re (\lambda) = 0$. This either implies a single root $\lambda = 0$, or a complex-conjugate root pair $\lambda = \pm {\rm i} \omega$. The case of a single real characteristic exponent $\lambda = 0$ is associated with `static' loss of stability. With this substitution we obtain
\begin{equation}\label{eq:D}
    D(\lambda = 0) = k_1 k_2 k_3 \sigma_4 + k_1 k_3 \sigma_4 + k_2 k_3 \sigma_4.
\end{equation}
One can find boundaries corresponding to pure complex characteristic exponents numerically, using the multi-dimensional bisection method \cite{MDBM}, {an extension of the bisection root-finding algorithm for multivariate implicit equation-systems}.

With the physical stiffness parameters $k_1, k_2, k_3$ being positive this would only be the case if the numerical stiffness parameter was chosen to be zero. This case can be discarded as realistic parameters ($\sigma_4 > 0$) are of interest.

Loss of stability may also happen through exponentially increasing oscillations if a complex conjugate pair of characteristic exponents crosses the imaginary axis ($\lambda = \pm {\rm i} \omega$). Such boundaries can be found by solving the transcendental equation $D({\rm i} \omega) = 0.$

Considering the physical substructure with the parameters given in Table~\ref{tab:id_par} and selecting the parameters $\gamma_4 = 3.7421$ Ns/m and $\sigma_4 = 45052$ N/m for the numerical damping and stiffness, boundaries were searched for in the parameter domain defined by the limits $\tau \in [5\times10^{-5}\,\text{s}, 30\times10^{-5}\,\text{s}]$, $p = \mu_3/m_3 \in [0.01,1]$, $\omega \in [10^{-4},550]$. These boundaries are shown in red in the stability charts in Fig.~\ref{fig:stabcharts}. 

To demonstrate that the results on the effect of imperfect interface control are  transferable to other hybrid testing scenarios, we present the stabilisability chart for the case when the numerical substructure has an extra DoF, thus simulating a five-storey building. We abstain from expressing the governing equations for this case; nevertheless, they are straightforward to derive based on those of the four-storey building. The stabilisability chart for the five-storey structure with the parameters $\sigma_5 = \sigma_4$, $\gamma_5 = \gamma_4$, $\mu_5 = m_1$ describing the additional top storey is presented in the right panel of Fig.~\ref{fig:stabcharts}. The structure of the stabilisability charts is very similar for the four- and five-storey structures suggesting that the instability is intrinsic to the imperfect interface control and similar effects are to be expected for other hybrid testing set-ups with a full-feedback control.

\begin{table}[!h]
    \caption{Identified parameters of the physical structure}
    \label{tab:id_par}
    \begin{tabular}{ l c c c }
        \hline
         & $m_r$ & $c_r$ & $k_r$ \\
        \hline
        1st storey ($r=1$) & 5.37 kg & 4.4062 Ns/m & 53048 N/m \\
        2nd storey ($r=2$) & 5.37 kg & 3.7979 Ns/m & 45724 N/m \\
        3rd storey ($r=3$) & 5.37 kg & 3.7421 Ns/m & 45052 N/m \\
        \hline
    \end{tabular}
    \vspace*{-4pt}
\end{table}

{Knowing the boundaries from Eq.~\eqref{eq:D} does not reveal the stability of the system in the resulting parameter-domains. Therefore, the different parameter-domains are labelled as stabilisable/unstable by semi-discretisation \cite{Insperger}} {a method of generating a set of (discretised) equations from the underlying infinite dimensional delay equation} using Chebyshev polynomials \cite{BredaDDEsemid,Trefethen2000}. Thus, the infinite dimensional time-delay system is approximated by a large-dimensional system of ordinary differential equations. This system can be given in the form $\dot{\mathbf{x}} = \mathbf{Ax}$, and the stability of the trivial equilibrium can be investigated by calculating the largest real-part of eigenvalues of the coefficient-matrix $\mathbf{A}$. {We say that the hybrid system is stabilisable if the trivial equilibrium in the semi-discrestised model is stable. This is due to the fact that in the real-life experiment, stability is also subject to using a suitable control algorithm to impose feedback from the numerical substructure on the physical substructure}. The unstable domains in Fig.~\ref{fig:stabcharts} are white whereas the stabilisable ones are indicated by green shading. The boundaries identified with the two different methods show a good agreement as there are only minor differences between the stabilisable domains identified with semi-discretisation and the curves obtained directly from the characteristic function $D(\lambda)$, which can be accounted for the semi-discretised system being an approximation of the time-delay model.

It is worth mentioning that both numerical methods have limitations in obtaining the stability boundaries at small time delays ($\tau \rightarrow 0$). This is due to the fact, that the semi-discretised approximation of the time-delay system has a singularity at $\tau = 0$. (One may overcome this shortcoming of the semi-discretisation by using linear multi-step methods to calculate the rightmost eigenvalue of the delay differential equation \cite{DDE_LMS}.) Another limitation of the numerical root-finding algorithms employed is that it requires a pre-defined domain in which solutions are searched for which makes it infeasible to find the stability boundaries from the characteristic equation, since solutions near $\tau = 0$ involve large frequencies ($\omega \rightarrow \infty$). Instead, for this case, we carry out an asymptotic estimation of the boundaries.

\section{Methodology -- Hybrid-testing without direct feedback}

As shown in Section 2, most of the difficulties in real-time dynamic substructuring of the 4-storey building arise in the feedback from the numerical to the physical substructure due to the actuator dynamics. Thus, we address these challenges by adopting a measurement strategy that breaks the feedback loop and replaces it with an iteration ensuring the synchronous motion of the physical and the numerical sides of the interface. This concept is presented in panel (b) of Fig.~\ref{fig:substr}.

The idea of using an iteration is adopted from the control-based continuation technique \cite{Sieber2008} which was originally developed to trace both stable and unstable periodic solutions in nonlinear experiments employing stabilising and non-invasive control to the observed structure{; the non-invasiveness of the control ensures that the steady-state behaviour of the system is identical to that of the open-loop system, changing the stability of steady-states but preserving their location in state-space}. In control-based continuation, the control target and the steady state system response is considered by their truncated Fourier series, allowing for the construction of an algebraic root-finding problem. Then an iteration is used to find the appropriate Fourier coefficients of the control target ensuring the non-invasiveness of the control, i.e., that the captured steady-state response of the controlled system is the same as the one of the open-loop system.

A similar strategy is proposed for the dynamic substructuring experiment. We take the iterative approach to find the steady-state response of the physical-numerical hybrid system to periodic forcing, i.e. its frequency response function while we disregard the transient response. Since, as a result, both the forcing and the steady-state system response will be periodic, the Fourier coefficients can be used to consider an algebraic zero-problem as in control-based continuation. This potentially allows for high-fidelity steady-state testing which is ideal in engineering applications where the FRF is more often of interest than the transients, as in product design or prototype testing, avoiding resonance is usually a fundamental requirement.

\newtext{We note that the iterative method described does not depend on the linearity of the system; the use of a Newton iteration allows nonlinear system interactions to be handled, as demonstrated in Section~\ref{sec:Results}.}

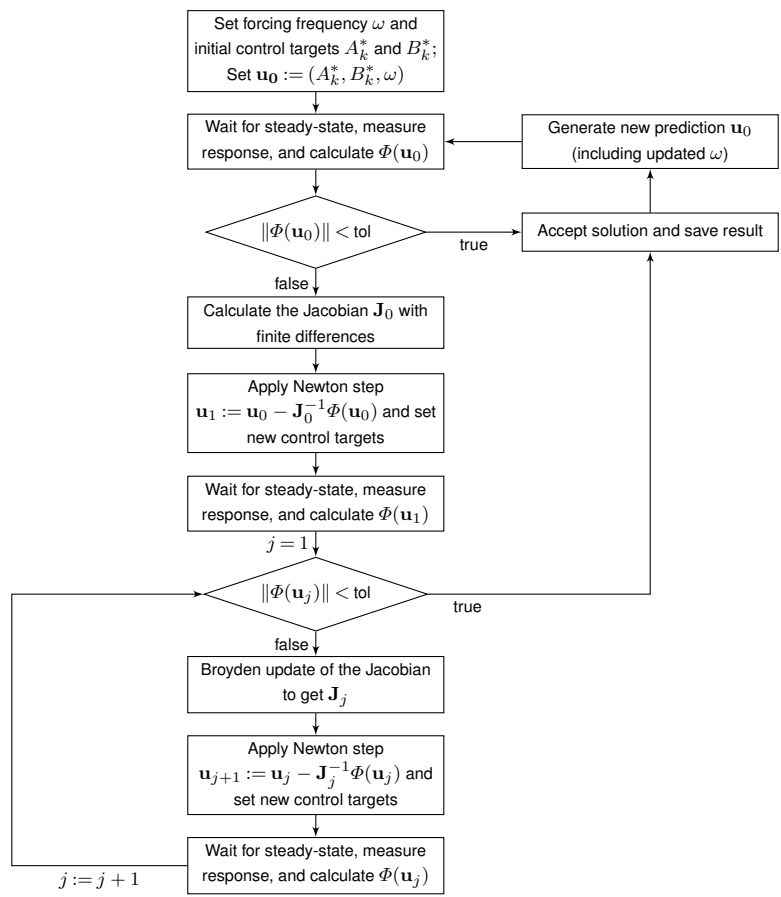
\begin{figure}[!h]
    \centering
    \begin{tikzpicture}[node distance = 15mm, font=\sffamily, scale=1, every node/.style={scale=0.8}]
        \node [process] (A) {\parbox{40mm}{\centering Set forcing frequency $\omega$ and initial control targets $A_k^*$ and $B_k^*;$ Set $\mathbf{u_0} := (A_k^*, B_k^*, \omega)$}};
        \node [process, below of=A] (B) {\parbox{40mm}{\centering Wait for steady-state, measure response, and calculate $\Phi(\mathbf{u}_0)$}};
        \node [decision, below of=B] (C) {$\lVert\Phi(\mathbf{u}_0)\rVert < \text{tol}$};
        \node [process, below of=C] (D) {\parbox{40mm}{\centering Calculate the Jacobian $\mathbf{J}_0$ with finite differences}};
        \node [process, below of=D] (E) {\parbox{40mm}{\centering Apply Newton step $\mathbf{u}_1 := \mathbf{u}_0 - \mathbf{J}_0^{-1}\Phi(\mathbf{u}_0)$ and set new control targets}};
        \node [process, below of=E] (F) {\parbox{40mm}{\centering Wait for steady-state, measure response, and calculate $\Phi(\mathbf{u}_1)$}};
        \node [decision, below of=F] (G) {$\lVert\Phi(\mathbf{u}_j)\rVert < \text{tol}$};
        \node [process, below of=G] (H) {\parbox{40mm}{\centering Broyden update of the Jacobian to get $\mathbf{J}_j$}};
        \node [process, below of=H] (I) {\parbox{40mm}{\centering Apply Newton step $\mathbf{u}_{j+1} := \mathbf{u}_j - \mathbf{J}_j^{-1}\Phi(\mathbf{u}_j)$ and set new control targets}};
        \node [process, below of=I] (J) {\parbox{40mm}{\centering Wait for steady-state, measure response, and calculate $\Phi(\mathbf{u}_j)$}};
        \node [process, right of=C, xshift=40mm] (K) {\parbox{40mm}{\centering Accept solution and save result}};
        \node [process, right of=B, xshift=40mm] (L) {\parbox{40mm}{\centering Generate new prediction $\mathbf{u}_0$ (including updated $\omega$)}};
        \path [connector] (A) -- (B);
        \path [connector] (B) -- (C);
        \path [connector] (C) -- node[anchor=east] {false} (D);
        \path [connector] (D) -- (E);
        \path [connector] (E) -- (F);
        \path [connector] (F) -- node[anchor=east] {$j=1$} (G);
        \path [connector] (G) -- node[anchor=east] {false} (H);
        \path [connector] (H) -- (I);
        \path [connector] (I) -- (J);
        \path [connector] (J) -- node[anchor=north] {$j := j + 1$} ++(-40mm, 0) |- (G);
        \path [connector] (C) -- node[anchor=north] {true} (K);
        \path [connector] (K) -- (L);
        \path [connector] (L) -- (B);
        \path [connector] (G) -| node [anchor=north, xshift=-30mm] {true} (K);
    \end{tikzpicture}

    \caption{Flowchart of the iterative substructuring technique.}
    \label{fig:flowchart}
\end{figure}

\subsection{Iterative root-finding and continuation}\label{ssec:rootf}

To match the physical and numerical displacements $x_3(t)$ and $\xi_3(t)$, which correspond to the interface storey of the assembly, we employ an iterative approach in the substructuring experiment. We consider a scenario when the hybrid structure is subject to periodic forcing (applied on the 4\textsuperscript{th} storey). As such, we expect a periodic response in steady state. This allows us to consider the displacements at the physical and numerical sides of the interface by means of their truncated Fourier series
\begin{equation}
    x_3(t) \approx A_0 + \sum_{k=1}^N A_k \cos (k \omega t) + \sum_{k=1}^N B_k \sin (k \omega t),
\end{equation}

\begin{equation}
    \xi_3(t) \approx \alpha_0 + \sum_{k=1}^N \alpha_k \cos (k \omega t) + \sum_{k=1}^N \beta_k \sin (k \omega t),
\end{equation}

If the tested structure is linear, then the higher harmonics ($k \geq 2$) in the steady-state response are small and practically negligible {(inevitably, non-idealities in the experiment generate some higher harmonics)}. The static components $A_0$ and $\alpha_0$ of the displacement signals are also discarded from the iteration. In case of the physical substructure, we assume that if the static part of the forcing from the shaker is not changing, the physical substructure oscillates around the static equilibrium whereas this is ensured by control in the numerical substructure.

In principle, the FRF of the hybrid system could be obtained by finding the appropriate phase and amplitude of the voltage signal sent to the actuator, resulting in the synchronous motion of the physical and numerical sides of the interface. However, since the physical structure is lightly damped, it is beneficial to apply a proportional derivative (PD) displacement control to the physical side of the interface and treat the control target $x^*(t)$ as an unknown to be determined by the iterative scheme. As such, introducing further damping to the structure, the transients are reduced, steady-state is reached faster, and the iteration can be performed with shorter waiting times. \newtext{The output of the controller $U_{\rm shaker}$ is used as a demand signal for a power amplifier, which drives the shaker. The power amplifier operates in current mode, and so the demand signal is (approximately) proportional to the force applied. It is not necessary to know the precise relationship between the demand signal and the force; the iteration is able to match the physical and numerical components at the interface regardless. However, the applied force $F_{\rm LC}$ is measured through a load cell and is used in the input to the numerical model.}

The demand signal sent to the actuator power amplifier is given by the PD control law
\begin{equation}
    U_{\rm shaker} = k_{\rm p} (x^*(t) - x_3(t)) + k_{\rm d} (\dot x^*(t) - \dot x_3(t))),
\end{equation}
where {the control target $x^*(t)$ is an arbitrary periodic function described by its Fourier series}
\begin{equation}
    x^* (t) = \sum_{k=1}^N A_{k}^* \cos (k \omega t) + \sum_{k=1}^N B_{k}^* \sin (k \omega t).
\end{equation}
{The constants $A_k^*$ and $B_k^*$ are determined as part of the iterative scheme, as described below} \newtext{where the index $k$ refers to the number of the harmonic component.} 

With these considerations, the aim of the iterative algorithm is to eliminate the error between the displacements of the physical and numerical 3\textsuperscript{rd} storeys such that

\begin{equation}
    \epsilon(t) = x_3(t) - \xi_3(t).
\end{equation}

Representing the steady-state solutions by their Fourier coefficients allows us to transform this into an algebraic problem where the error vector $\boldmath{\Phi}$ is defined by the difference of the corresponding Fourier coefficients
\begin{equation}\label{eq:fcoeff}
    \Phi (A_k^*, B_k^*) = 
    \left(\begin{matrix}
        A_k - \alpha_k\\
        B_k - \beta_k
    \end{matrix}
    \right)
    = \left( \begin{matrix}
        A_1 - \alpha_1\\
        \vdots\\
        A_N - \alpha_N\\
        B_1 - \beta_1\\
        \vdots\\
        B_N - \beta_N
    \end{matrix}\right).
\end{equation}
{Here, $A_k^* = (A_1^*, \ldots, A_N^*)^\text{T}$ and $B_k^* = (B_1^*, \ldots, B_N^*)^\text{T}$ are the Fourier coefficients of the control target. The Fourier coefficients corresponding to the physical displacement $(A_k, B_k)$ and the numerical displacement $(\alpha_k, \beta_k)$ are both functions of $(A_k^*, B_k^*)$, defined implicitly through the behaviour of the physical substructure. Each time a new control target is set, the hybrid experiment equilibrates and new values for $(A_k, B_k)$ and $(\alpha_k, \beta_k)$ are measured.} This error vector $\Phi$ is already suitable for ensuring the synchronous motion of the physical and the numerical 3\textsuperscript{rd} storeys for a constant set of system parameters. 

If the frequency response of the structure is a well-defined function (i.e. single valued) with respect to a system parameter of interest (e.g., the forcing frequency), the parameter domain can be swept to determine the overall behaviour. However, in the presence of weak nonlinearity hysteresis loops may form leading to a multi-valued frequency response. \newtext{In this case, pseudo-arclength continuation can be applied directly to the hybrid experiment, see for example~\cite{Barton2017}, in which the solution curve is re-parameterised by the (pseudo) arclength rather than using the forcing frequency. The addition of the pseudo-arclength condition means that the forcing frequency is treated as another unknown that is solved for as part of the continuation process, and so during solving the experiment will be run for varying values of the forcing frequency. This enables solution curves to be tracked through the fold points, which give rise to multi-valued responses.}

\newtext{Since pseudo-arclength continuation is not required for the implementation of the iterative hybrid testing approach, we leave a detailed explanation of the method to excellent textbooks such as~\cite{Seydel}. Alternatively, there are off-the-shelf software packages such as \textsc{CoCo}~\citep{Dankowicz2013} and BifurcationKit.jl~\citep{Veltz2020} that implement pseudo-arclength continuation already.}

The pseudo-arclength condition is appended to~\eqref{eq:fcoeff} to yield the extended system 
\begin{equation}\label{eq:fcoeff_pseud}
    \Phi_{\rm ext} (\mathbf{u}) = \Phi_{\rm ext} (A_k^*, B_k^*, \omega) = \left( \begin{matrix}
        A_k - \alpha_k \\
        B_k - \beta_k \\
        \mathbf{u}_{\rm secant} \cdot (\mathbf{u}_{\rm pred} - \mathbf{u}) 
    \end{matrix}\right),
\end{equation}
where $\omega$ is a system parameter and the input vector is given by $\mathbf{u} = (A_k^*,B_k^*,\omega)^\text{T}$. The secant $\mathbf{u}_{\rm secant} = \mathbf{u}_{-1} - \mathbf{u}_{-2}$, is obtained by using the two previous solutions, is used to provide the predicted direction for the continuation of the solution branch. \newtext{To initialise this method, we find the first two points on the branch without the pseudo-arclength condition solving \eqref{eq:fcoeff} only.} As such, the prediction of the next solution is calculated as $\mathbf{u}_{\rm pred} = \mathbf{u}_{-1} \Delta s\, \mathbf{u}_{\rm secant}$, where the constant $\Delta s$ is used for step-size control. \newtext{In our implementation we use $\Delta s$ = 1.1 if the previous iteration converged to a new point on the solution curve. If, however, the root finding iteration does not converge, the step size is set to $\Delta s = 0.5$ and the step is repeated.}

Starting from a suitable initial condition $\mathbf{u}_0$, the solution of the extended zero-problem $\Phi_{\rm ext}(\mathbf{u}) = \mathbf{0}$ is found by a Newton-like iteration
\begin{equation}\label{eq:nexstep}
    \mathbf{u}_{j+1} = \mathbf{u}_{j} - \mathbf{J}^{-1}_j \Phi(\mathbf{u}_j).
\end{equation}
In the original Newton method for nonlinear root-finding, $\mathbf{J}_j$ stands for the Jacobian matrix. However, this is unavailable in experiments; therefore, a suitable proxy is used in the algorithm. At the first iteration step, $\mathbf{J}_1$ is obtained by finite differences, which equals the same number of evaluations of the error vector $\Phi$ (obtaining steady-state solutions for different inputs) as the problem dimension. Then, in the following steps ($j \geq 2$) a Broyden update \cite{Broyden} is used to the Jacobian which requires one evaluation at each step.

\begin{equation}
    \mathbf{J}_{j+1} = \mathbf{J}_{j} + \frac{(\Phi(\mathbf{u}_j)-\Phi(\mathbf{u}_{j-1})) - \mathbf{J}_{j} (\mathbf{u}_{j} - \mathbf{u}_{j-1})}{(\mathbf{u}_{j} - \mathbf{u}_{j-1})^2}(\mathbf{u}_{j} - \mathbf{u}_{j-1})^T.
\end{equation}

Since Newton-like root-finding methods are prone to overshooting (i.e., not converging due to taking too large steps) the algorithm is combined with line-search modifying the iteration step in Eq.~\eqref{eq:nexstep} to
\begin{equation}
    \mathbf{u}_{j+1} = \mathbf{u}_{j} - \nint \mathbf{J}^{-1}_j \Phi(\mathbf{u}_j).
\end{equation}
Initially, the next candidate solution is searched for with $\nint = 1$. If then $\lVert \mathbf{u}_{j+1} \rVert_2 < \lVert \mathbf{u}_{j} \rVert_2$ this solution is accepted, if not, we retry with $\nint := \nint/2$ until a smaller error norm is achieved than in the previous step or the pre-set minimum value of the damping parameter $\nint$ is reached. 

The iteration is either terminated upon reaching the given tolerances or exceeding the maximum number of iterations, in which case the solution is discarded. In this case, the continuation may progress by considering a smaller step-size from the previous point. The flowchart of the iterative method can be seen in Fig.~\ref{fig:flowchart}.

\section{Experimental investigation}

\subsection{Hybrid testing setup}

In the dynamic substructuring experiment, we investigate a physical-numerical hybrid model of a bench-top-sized four-storey building. As shown in Figure~\ref{fig:substr_with_exp}. The physical substructure comprises the two bottom storeys and the legs of the physical structure whereas the third storey, attached to the shaker, serves as the interface between the physical and numerical substructures. The displacements of the three storeys of the physical substructure are measured by laser sensors. A piezoelectric force transducer is used to measure the force input from the shaker to the physical part which is fed to the real-time simulation of the numerical substructure. The acceleration of the three physical storeys are captured by accelerometers the signals of which can be used for compensating the inertial forces. This configuration allows for testing different mass-ratios of the physical and numerical parts of the third storey of the true assembly.

The numerical substructure comprises the third and fourth storeys of the building. The real-time simulation of the numerical model is realised by the real-time controller (RTC) unit which is also used for data acquisition. The RTC device is also responsible for controlling the shaker providing the control force to the physical system.

The RTC is connected to a PC which is used for data processing and determining the control target and parameters.

As shown in Fig.~\ref{fig:substr_with_exp}, in the hybrid testing setup, the connection of the physical and numerical substructures is realised by feeding the measured force at the interface to the numerical model, whereas the simulated displacement of the numerical third storey is imposed on the physical substructure by the shaker.

\subsection{Simulation of the numerical substructure}

The numerical substructure, as indicated in Fig.~\ref{fig:substr_with_exp} , comprises a part of the 3$^{\rm rd}$ and the 4$^{\rm th}$ storeys of the building and the columns in between. This part of the structure is simulated by a 2 DoF numerical model which is subject to three external forces; the interface force $F_{\rm LC}(t)$, measured by the load cell in real-time, the control-force $F_{\rm ctrl}$ acting on the 3$^{\rm rd}$ storey and, in the present measurement scenario, a periodic forcing on the 4$^{\rm th}$ storey $F_{\rm forc} = \Phi \cos(\omega t)$. The control-force $F_{\rm ctrl}$ is employed on the numerical structure to eliminate any static error in displacement due to the numerical structure being ungrounded and is designed such that it is non-invasive, i.e. it does not affect the steady state response of the hybrid structure.

Furthermore, to carry out hybrid testing for different mass ratios $p=\mu_3/m_3$, we use the measured acceleration signal to compensate for a part of the inertial force coming from the physical side of the interface $F_{\rm comp} = m_{\rm 3 comp} \ddot x_3$. 

\newtext{We also apply proportional feedback control to the numerical part of the structure to ensure the third (interface) storey oscillates around $\xi_3 = 0$. This is necessary since the numerical substructure is not grounded. Without this grounding, the static error of the piezoelectric force and acceleration sensors as well as the error in the numerical integration would lead to an error in displacement. The control-force applied to the numerical part of the third storey is determined by the PD control law}
\newtext{
\begin{equation}
    F_{\rm ctrl} = \hat k_{\rm p} (x_3(t) - \xi_3(t)) + \hat k_{\rm d} (\dot x_3(t) - \dot \xi_3(t)).
\end{equation}
}
\newtext{This definition ensures that the control applied to the numerical substructure is non-invasive, i.e., if the iteration ensures that the physical and the numerical displacements of the interface storey are equal in steady state, the control force $F_{\rm ctrl}$ vanishes and the response of the hybrid structure is equal to the one of the true assembly.}

Thus, the governing equations of the numerical model can be expressed as
\begin{equation}
    \mu_3 \ddot \xi_3(t) = F_{\rm num} + \gamma_4 (\dot \xi_4(t) - \dot \xi_3(t)) + \sigma_4 (\xi_4(t) - \xi_3(t))
\end{equation}
and
\begin{equation}
    \mu_4 \ddot \xi_4(t) = \gamma_4 (\dot \xi_3(t) - \dot \xi_4(t)) + \sigma_4 (\xi_3(t) - \xi_4(t))+F_{\rm forc}(t),
\end{equation}
where the external forces acting on the numerical 3\textsuperscript{rd} storey are given by
\begin{equation}
    F_{\rm num}(t) = F_{\rm comp}(t)- F_{\rm LC}(t) + F_{\rm ctrl}(t).
\end{equation}

Introducing new state variables $z_1 = \dot \xi_3$, $z_2 = \dot \xi_4$ ,$z_3 = \xi_3$, $z_4 = \xi_4$, the governing equations of the numerical substructure can be re-written in the matrix form $\dot{\mathbf{z}} = \mathbf{A} \mathbf{z} + \mathbf{b}$

\begin{equation}
    \left( \begin{matrix} \dot z_1 \\ \dot z_2 \\ \dot z_3 \\ \dot z_4 \\ \end{matrix} \right) = \left( \begin{matrix}
        -\frac{\gamma_4}{\mu_3} & \frac{\gamma_4}{\mu_3} & -\frac{\sigma_4}{\mu_3} & \frac{\sigma_4}{\mu_3} \\
        \frac{\gamma_4}{\mu_4} & -\frac{\gamma_4}{\mu_4} & \frac{\sigma_4}{\mu_4} & -\frac{\sigma_4}{\mu_4} \\
        1 & 0 & 0 & 0 \\
        0 & 1 & 0 & 0 \\
    \end{matrix}\right) \left( \begin{matrix} z_1 \\ z_2 \\ z_3 \\ z_4 \\ \end{matrix} \right) + \left( \begin{matrix} \frac{F_{\rm LC}}{\mu_3} \\ \frac{F_{\rm forc}}{\mu_4} \\ 0 \\ 0 \\ \end{matrix} \right).
\end{equation}

This system is simulated using the implicit Euler scheme and so the subsequent values of the state variables are given by 
\begin{equation}
    \mathbf{z}_{i+1} = \mathbf{z}_{i} + h \left( \mathbf{A} \mathbf{z}_{i+1} + b_{i+1}\right),
\end{equation}
with $h$ referring to the simulation time-step.
Since the numerical substructure is linear, one can solve this equation for $\mathbf{z}_{i+1}$ obtaining
\begin{equation}
    \mathbf{z}_{i+1} = \left( \mathbf{I} - h \mathbf{A} \right)^{-1} \left( \mathbf{z}_{i} - h \mathbf{b}_{i} \right).
\end{equation}
\newtext{While for this analysis, it would be satisfactory to compute the frequency response, the apporach of performing a full time-domain simulation in the numerical substructure more versatile, as it could be extended to study nonlinear structures or non-periodic inputs.}

The simulation of the numerical substructure is realised using the parameters listed in Table~\ref{tab:virt_par}. 

\begin{table}[!h]
    \caption{Parameters of the numerical structure}
    \label{tab:virt_par}
    \begin{tabular}{ l l c }
        \hline
        $\mu_4$ & 4\textsuperscript{th} storey mass & 10.74 kg\\
        $\gamma_4$ & Numerical leg damping & 3.7421 Ns/m \\
        $\sigma_4$ & Numerical leg stiffness & 45052 N/m \\
        $\Phi$ & Forcing amplitude & 1 N \\
        $h$ & Simulation timestep & $10^{-4}$ s \\
        $k_p$ & Physical proportional control gain & 100 V/m\\
        $k_d$ & Physical derivative control gain & 10 Vs/m\\
        $\hat k_p$ & Numerical proportional control gain & 50000 N/m\\
        $\hat k_d$ & Numerical derivative control gain & 200 Ns/m\\
        \hline
    \end{tabular}
    \vspace*{-4pt}
    \end{table}

\section{Results}
\label{sec:Results}

\subsection{Frequency response of the hybrid system}

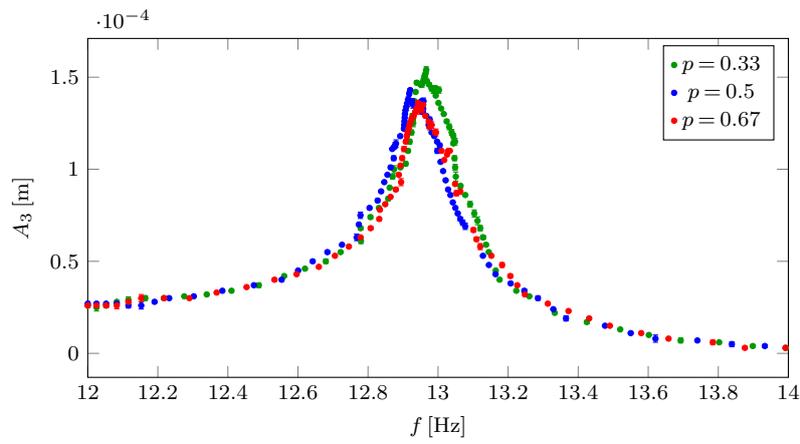
\begin{figure}[]
    \centering
    \begin{tikzpicture}
\begin{axis}[
legend style = {font=\mtcaption},
xlabel={$f \; {\rm [Hz]}$},
ylabel={$A_3 \; {\rm [m]}$},
mark size = 1pt,
xmin = 12,
xmax = 14,
width=0.8\textwidth,
height=9/16*0.8\textwidth
]
\addplot [black!40!green, only marks, error bars/.cd,
y dir=both,y explicit relative] table[x=x,y=y,y error=error] {figs/freq_sweep_0.33.dat};
\addlegendentry{$p=0.33$}
\addplot [blue, only marks, error bars/.cd,
y dir=both,y explicit relative] table[x=x,y=y,y error=error] {figs/freq_sweep_0.5.dat};
\addlegendentry{$p=0.5$}
\addplot [red, only marks, error bars/.cd,
y dir=both,y explicit relative] table[x=x,y=y,y error=error] {figs/freq_sweep_0.66.dat};
\addlegendentry{$p=0.67$}
\end{axis}
\end{tikzpicture}
    \caption{Frequency response of the physical-virtual hybrid 4-storey building with different mass-ratios at the interface. \newtext{The frequency sweep is performed around the second bending mode of the structure.} Error bars show the magnified (2x) error between the steady-state amplitudes of the fundamental harmonic components of the physical and virtual third-storey displacements.}
    \label{fig:4_storey_lin_resp}
\end{figure}

Using the continuation algorithm, described in Section 3.~\ref{ssec:rootf}, the frequency response of the 4-storey building was captured for periodic forcing applied to the 4\textsuperscript{th} storey. {The second bending mode of structure was considered, since the structure is effectively linear in this configuration, whereas the first bending mode showed noticeable nonlinearity. As such,} only the fundamental harmonic components of the physical and numerical displacements were matched by the iterative algorithm (using $k = 1$ in Eq. \eqref{eq:fcoeff} and \eqref{eq:fcoeff_pseud}). Using different masses in the numerical 4\textsuperscript{th} storey and compensating for the inertial forces accordingly, we tested three different mass ratios $p=\mu_3/m_3$ of the numerical and physical sides of the interface: Case 1: p = 0.67 ($\mu_3 = m_{\rm 3 comp} = 3.58$ kg), Case 2: p = 0.5 ($\mu_3 = m_{\rm 3 comp} = 2.69$ kg), p = 0.33 ($\mu_3 = m_{\rm 3 comp} = 1.79$ kg). Meanwhile, the total mass of the 3\textsuperscript{rd} storey was kept constant $m_3 + \mu_3 = 5.37$ kg. As such, the three cases considered here can be interpreted as three different interface position, with the mass of the 3\textsuperscript{rd} storey split between the physical and numerical substructures at different heights.

We draw attention to the fact that, as presented in Section 4, a necessary condition of stabilisability for closed-loop feedback control is $p>1/2$. As such, if a traditional closed-loop feedback control was used for substructuring, Case 1 would be stabilisable, Case 2 corresponds to the theoretical lower limit for stabilisable mass ratios whereas Case 3 would be unconditionally unstable. In contrast, with the iterative approach, it was possible to capture the FRF for all cases irrespective of the interface position. Moreover, the good agreement between the cases with different mass ratios indicate that the algorithm results in minimal interface error and the results are representative of the behaviour of the true assembly.

Figure~\ref{fig:4_storey_lin_resp} shows the FRFs indicating the displacement amplitude for the 3\textsuperscript{rd}, interface storey. The iterative technique is capable of finding the solution within $0.1\%$ or $10^{-4}$\,m tolerance (whichever is the higher) for the Fourier coefficients of fundamental harmonic component of the solutions, resulting in a very accurate match between the physical and numerical displacements. {If error bars were used for Fig.~\ref{fig:4_storey_lin_resp}, as in case of the nonlinear frequency response curve in Fig.~\ref{fig:nl_freq_resp}, they would be smaller or similar in size to the markers in the diagram}. Note that the experiment has less overall accuracy in terms of {independent} repeatability of the FRF. The overall errors however are intrinsic to the experiment {(such as daily temperature changes and imperfect vibration isolation within the lab)} and not introduced by the iterative scheme.

Identifying a valid data point in the FRF while performing the continuation typically requires 0-5 iterations for any mass-ratio, though there were occasional outliers which can be explained by the larger noise level the experiment is subject to at larger amplitudes. It is also worth noting that there is no significant difference in the performance of the iterative algorithm for the three different mass ratios, indicating that the iterative substructuring scheme performs evenly in the three cases. This is in contrast with the direct feedback-based approach, where e.g. the case $p=0.33$ would be theoretically impossible due to the actuator delay. 

\subsection{Frequency response of a structure with nonlinear restoring force}

\begin{figure}[!t]
    \centering\includegraphics{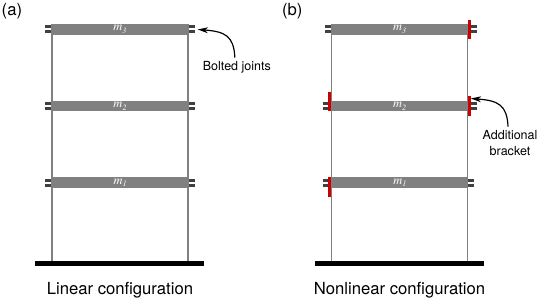}
    \caption{A schematic diagram of the physical substructure in both (a) linear and (b) nonlinear configurations. In the nonlinear configuration, the additional brackets bolted onto the masses provide a bi-linear stiffness.}
    \label{fig:nonlinear-schematic}
\end{figure}

\begin{figure}[!t]
    \centering\includegraphics[width=\textwidth]{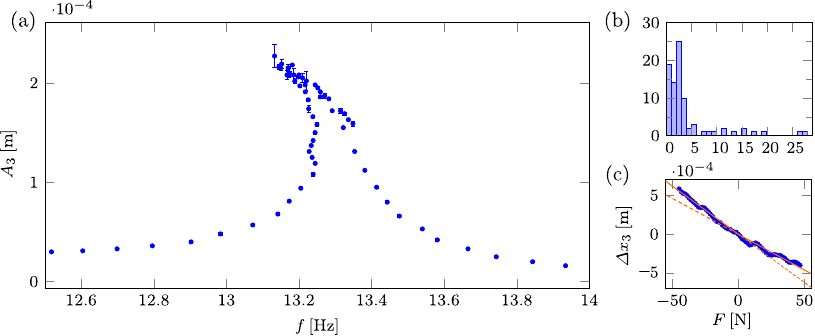}
    \caption{Panel (a): Frequency response of the physical-numerical hybrid 4-storey building with nonlinear restoring force (softening characteristics). Error bars show the magnified (2x) error between the amplitudes of the fundamental harmonic components of the physical
    and numerical third-storey displacements. Panel (b): Histogram indicating the number of iterations required to capture each recorded point of the FRF. Panel (c): Measured static force response (blue markers) and fitted bi-linear characteristics (continuous orange lines) with both linear sections extended along the axis (dashed orange lines).}
    \label{fig:nl_freq_resp}
\end{figure}

To further demonstrate the versatility of the iterative approach in substructuring, we performed hybrid testing with a modified physical substructure producing a nonlinear restoring force. To achieve this, stiff plates were added to the legs of the physical 3-storey structure to limit their effective length in one bending direction (see Fig.~\ref{fig:nonlinear-schematic}(b)). This resulted in nonlinear stiffness as indicated in panel (c) of Fig.~\ref{fig:nl_freq_resp}, showing the displacement difference $\Delta x_3 = x_3-x_2$ between the second and third storeys for a range of static forcing. The measured nonlinear response is approximated by bilinear force characteristics which is consistent with a different effective length when the relative displacement is in one direction compared to the other direction. In order to trigger a nonlinear behaviour, a static force of 26.68 N was subjected to the structure by the shaker, leading to softening  behaviour at larger amplitudes.

Panel (a) of Fig.~\ref{fig:nl_freq_resp} shows the FRF of the nonlinear structure with bars indicating magnified (2x) amplitude errors between the physical and numerical sides of the interface. Observe that the diagram includes a parameter domain where multiple solutions correspond to the same frequency. While the stability of the captured solutions is not investigated in the experiment, nonlinear system theory suggests that the high and low amplitude solutions are stable, whereas the solution-branch in the middle is unstable. Extracting information of such unstable limit cycles in nonlinear experiments is non-trivial as they cannot be found without control. Since the iterative substructuring method is derived from CBC it has the advantage over the direct feedback-based method that it is capable of finding even unstable limit cycles.

\begin{figure}[!h]
    \centering\includegraphics[]{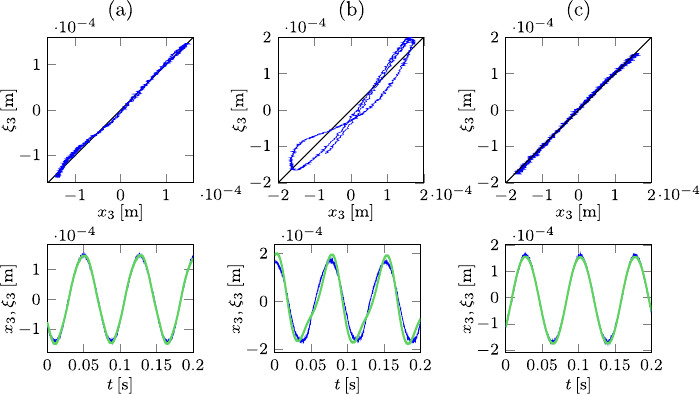}
    \caption{Top row: Phase plots of the physical and numerical interface displacements. The black line indicates the ideal case ($x_3(t) = \xi_3(t)$). Bottom row: time profiles of the physical (blue) and numerical (green) interface displacements. Column (a): linear physical substructure. Column (b): Bi-linear stiffness in the physical substructure, with only the fundamental harmonic components synchronised. Column (c): Bi-linear stiffness in the physical substructure, with synchronisation of the 2nd and 3rd harmonics \newtext{(i.e., 2nd and 3rd harmonics are included in \eqref{eq:fcoeff})}.}
    \label{fig:iteartions}
\end{figure}

In Fig.~\ref{fig:iteartions}, we compare the converged, steady-state solutions from the linear (column a) and the nonlinear (column b) case. In the two diagrams, we show two phase-plots of the physical and numerical displacements $x_3$ and $\xi_3$ in the top and bottom row of panels, respectively. In column (a), indicating a solution from the linear structure, only minor differences can be observed between the physical and numerical displacements, which are present mostly due to measurement and/or process noise or a weak effect of the higher harmonic components. Column (b) shows a case with more profound nonlinear behaviour. In this case, the fundamental harmonic components of the numerical and physical displacements are still matched but there is more prevalent difference in the higher (especially the second) harmonic as indicated by shape of the Lissajous plot \cite{Lissajous}. These differences can be eliminated by further iterations incorporating the higher harmonic components also, as indicated by column (c), captured with the nonlinear physical substructure and at a similar vibration amplitude as in case (b). However, matching the higher harmonics takes significantly longer than just the fundamental harmonic component of the solution. Thus, this was not used when capturing the nonlinear FRF.

Note that the realisation of the nonlinear restoring force in the presented experiment has some limitations. Since the bi-linear response originates from contact-nonlinearity, the experiment is subject to a significantly higher level of noise and uncertainty in the nonlinear configuration than in the linear case. This is reflected by the histogram in panel (b) of Fig.~\ref{fig:nl_freq_resp} showing that, while in most cases the solution was found in 5 iterations, there were a considerable number of instances when a higher number of iteration steps was required. This is due to the higher fluctuation in the captured response, which limited the use of the pseudo-arclength continuation resulting in slightly poorer coverage of the response curve. {As indicated by the larger error bars}, these effects were strongest around the resonance peak. Therefore, the nonlinear response was captured by two continuations one approaching the peak from lower frequencies while the other from higher frequencies. Moreover, the higher fluctuation in the Fourier coefficients prevented the inclusion of the higher harmonics to the iteration, as since a higher-dimensional Jacobian matrix is increasingly difficult to obtain. Note that these issues could be addressed by employing specialist techniques (see \cite{SCHILDER2015}) for performing CBC in noisy or impacting systems involving higher uncertainty. Nevertheless, this was deemed outside the scope of the present study.  

\section{Conclusion}
In this study, the dynamic substructuring of a bench-top-sized four-storey building is performed with the first two storeys fully physically realised, the third storey serving as the interface to the numerical substructure, and the top two storeys simulated. Using a time-delayed model to account for actuator dynamics, our study demonstrates that the hybrid experiment is very challenging and in some cases impossible to carry out with the classical approach using a full feedback-loop.

Instead, an iterative approach is proposed, breaking the feedback-loop between the numerical and the physical parts of the system, thus avoiding the related stability problems. We demonstrated that the iterative technique is capable of capturing the FRFs of the hybrid system with good accuracy without the risk of the feedback-related instability. The iterative approach also avoids the use of model-prediction techniques or other more sophisticated means of compensation for the actuator dynamics. Being a model-free technique also has the advantage that the iterative technique can handle various changes or different model configurations without major alternations to the algorithm. Thus, it is a versatile tool that may be employed when more traditional substructuring approaches fail. As the method is based on numerical continuation, it can be easily used to track the response of nonlinear functions, even unstable solutions which would be undetectable otherwise.

Nevertheless, certain scenarios are still difficult to handle, i.e. solutions or measurement setups involving a higher amount of noise, making it difficult to find the solution resulting in the synchronous motion of the physical and numerical sides of the interface by Newton-like iteration, as derivatives are especially sensitive to random perturbations. In further studies, these issues may be addressed by using surrogate models relying on sampled data from the experiment. Also, the effect of noise and model uncertainty can be addressed providing confidence intervals for the captured FRFs for the given noise level and uncertainty in the system. {Finally, this approach has natural extensions to systems with multiple degrees-of-freedom at the interface; initial experiments can be found in~\cite{Alessandra_Hybrid}.}

\vskip6pt

\enlargethispage{20pt}

\ethics{The authors declare that this research and the presentation of its results is in line with the ethical standards of the Royal Society, no human or animal subjects were used and there are no ethical concerns related to biosecurity involved.}

\dataccess{Data are available at the University of Bristol data repository, \url{data.bris}, at \url{https://doi.org/10.5523/bris.38u46si2df1xf28ugj0b7cm15z}.}

\aucontribute{All authors contributed to the algorithm, experiment and manuscript development. S. Beregi carried out the experiments and wrote the manuscript.}

\competing{There are no competing interests.}

\funding{This research has received funding from the Digital twins for improved dynamic design (EP/R006768/1) EPSRC, United Kingdom grant. The support of the EPSRC is greatly acknowledged.}


\vskip2pc

\bibliographystyle{plainnat}
\bibliography{Hybrid_testing}



\appendix{}

\section{Asymptotic calculation of the stability boundary as $\tau \rightarrow 0$}
\label{sec:asymptotics}

Although an analytical solution is not available for the stability boundaries, an asymptotic calculation can be carried out for the case of oscillatory loss of stability ($\lambda = {\rm i} \omega$) for $\omega \rightarrow \infty$ and $\tau \rightarrow 0$.

\newtext{This means that in this case, the terms with the coefficient $\lambda^2$ outweigh all other terms in \eqref{eq:char_fun} and the characteristic function $D$ tends to}
\begin{equation}
    D(\omega) \rightarrow {\rm det} \left( -\omega^2 \mathbf{M}_{\tau} e^{-{\rm i} \omega \tau} -\omega^2 \mathbf{M}_0 \right).
\end{equation}

Using Euler's formula for the complex exponential this can be expressed as
\begin{equation}
    {\rm det} \left( {\rm i} \omega^2 \mathbf{M}_{\tau} \sin (\omega \tau) - \omega^2 \left(\mathbf{M}_{\tau} \cos (\omega \tau) + \mathbf{M}_0 \right) \right) = 0.
\end{equation}
If $\omega \tau = (2k+1) \pi$ with $k \in \mathcal{Z}$ (although $\omega$ is very large while $\tau$ is very small, their product may be in leading order); then $\sin (\omega \tau) = 0$ and $\cos (\omega \tau) = -1$. Thus, the equation above yields

\begin{equation}
    {\rm det} \left( -\mathbf{M}_{\tau} +\mathbf{M}_{0} \right) = 0
\end{equation}
which leads to
\begin{equation}
    -\mu_3 + m_3 = 0.
\end{equation}
As indicated in Fig.~\ref{fig:stabcharts}, this analytical boundary is in good agreement with the numerical results, and it may also serve as a rule of thumb: substructuring with delayed feedback is infeasible if the mass in the physical part is larger than in the numerical part, i.e. if $p \leq 1/2$. Nevertheless, the numerical results indicate that the maximum time-delay allowed in the feedback loop is very limited (in the order of 0.1\,ms) even for the mass-ratios where substructuring is theoretically possible. This means that realising a dynamic hardware-in-the-loop test with a close feedback loop is very challenging and essentially has to rely on the compensation of the actuator dynamics.

\end{document}